\documentclass[11pt,letterpaper]{amsart}
\usepackage{amssymb}
\usepackage{enumitem}
\usepackage[colorlinks = true,
            linkcolor = red,
            urlcolor  = blue,
            citecolor = blue,
            anchorcolor = blue]{hyperref}
\newtheorem{theorem}{Theorem}[section]

\newtheorem{lemma}[theorem]{Lemma}
\newtheorem{corollary}[theorem]{Corollary}

\theoremstyle{definition}

\theoremstyle{remark}
\newtheorem{remark}[theorem]{Remark}

\newcommand{\Prob}{\mathbf{P}}

\numberwithin{equation}{section}
\usepackage{amsrefs}
\def\MR#1{\href{http://www.ams.org/mathscinet-getitem?mr=#1}{MR#1}}
\newcommand{\ARXIV}[1]{\href{https://arXiv.org/abs/#1}{arXiv:#1}}
\begin{document}

\title[Laws of large numbers for real roots of random polynomials]{A  strong law of large numbers for real roots of random polynomials}      
   \author{Yen Q. Do}
\address{Department of Mathematics, University of Virginia,  Charlottesville, Virginia 22904, USA}
\email{yqd3p@virginia.edu}

\keywords{Random polynomials, real roots,  strong law of large numbers, almost surely, convergence}
\subjclass[2020]{60G50, 60F05, 41A60}
\date{\today}
\begin{abstract} We consider random polynomials $p_n(x)=\xi_0+\xi_1+\dots+\xi_n x^n$ whose coefficients  are independent and identically distributed with zero mean,  unit variance, and  bounded $(2+\epsilon)^{th}$ moment (for some $\epsilon>0$),  also known as the Kac polynomials.   Let $N_n$ denote the number of real roots of  $p_n$.  In this paper, motivated by a question from Igor Pritsker,  we prove that almost surely the following convergence holds:
\begin{eqnarray*}
\lim_{n\to\infty} \frac{N_n([-1,1])}{\log n} &=& \frac 1 \pi.
\end{eqnarray*}
This convergence could be viewed as a local strong law for the real roots.  The main ingredient in the proof is a set of maximal inequalities  that reduces the proof to proving convergence along  lacunary subsequences, which in turn follows from a recent concentration estimate of Can--Nguyen.
\end{abstract} 
\maketitle


\section{Introduction}
The classical strong law of large numbers (SLLN) \cite{Dur} asserts that, for a sequence of  independent and identically distributed (iid) random variables $\eta_1,\dots,\eta_k,\dots $ with nonzero mean and finite first absolute moment ($\mathbb E|\eta_j|<\infty$),  almost surely we have
\begin{eqnarray*}
\lim_{k\to\infty} \frac{\eta_1+\dots+\eta_k}{\mathbb E [\eta_1 +\dots +\eta_k]} &=& 1.
\end{eqnarray*}

Our goal  is to investigate a similar law for the number of real roots for random polynomials 
\begin{eqnarray*}
p_n(x) &=& \xi_0+\xi_1 x+\dots + \xi_n x^n
\end{eqnarray*}
where  $\xi_0, \xi_1, \dots$ are iid with zero mean and unit variance ($p_n$ is also known as the Kac polynomials).  

Let $N_n$ be the number of real roots of $p_n$.  The distribution of the real roots of $p_n$ and in particular of $N_n$ have attracted the attention of mathematicians for several centuries, with seminal contribution of Kac \cite{K,K1}, Littlewood--Offord \cite{LO1,LO2, LO3},  Ibragimov--Maslova \cite{IM1,IM2,IM3,M,M1},  and Tao--Vu \cite{TV},  to name a few,  see also \cite{anw, ADL, AL,BCP, BS, BD, BP, CN, CP, C,Da, DPSZ,D1, D2, DM, DL, DLNNP, DNNP, DNN,DV, Dun1, Dun2, Dun3, EK, EO, F, G,GW,HZ, HKPV, KZ, Kost,KLN, LP, LP1, LPX, LPX1, NS1, NS2, N,NNV, NV1,NV2, Q,S, SGF,  STB, SM0, SM, St,W}. The expected value of $N_n$ is well-known: the classical work of Kac \cite{K} from 1940s first established  that 
\begin{eqnarray*}
\mathbb E N_n &=& \frac 2\pi \log n +o(\log n)
\end{eqnarray*} 
for random polynomials with standard Gaussian coefficients,  and this result has been extended to more general classes of coefficient distributions and with sharper bounds for the error term, for example assuming zero mean, unit variance, and finite $(2+\epsilon)^{th}$ moment we can show that the error term is $O(1)$  \cite{NNV}.   One may in fact divide $\mathbb R$ in to roughly $\log n$ disjoint regions, such that the average number of real roots inside each region is bounded, see e.g. \cite{TV,DNV}. Furthermore,  the numbers of roots inside these local regions essentially behave like weakly correlated random variables, and this intuition is implicit in many studies in this area, for example in Maslova's proofs \cite{M,M1} of the variance asymptotics $Var[N_n] =\frac 4 {\pi}(1-\frac 2{\pi})C\log n +o(\log n)$ and the convergence in distribution
\begin{eqnarray*}
\frac{N_n-\mathbb E N_n}{\sqrt{Var[N_n]}} &\to& N(0,1),
\end{eqnarray*}
also known as the central limit theorem for $N_n$ (see also \cite{NV1,DN} for recent extensions of these  celebrated results of Maslova to generalized Kac polynomials,  which include derivatives of the Kac polynomials and hyperbolic random polynomials).

The above discussion suggests that $N_n$ is very similar to a sum of weakly correlated random variables of the same typical order, and furthermore the number of terms grows to $\infty$ when $n\to\infty$. The classical SLLN  suggests that $N_n/\mathbb E[N_n]$ should converge almost surely to $1$, as $n\to\infty$. In fact, this conjecture was brought forward to the community  by Igor Pritsker during  the 2019 AIM workshop on ``Zeros of random polynomials" \cite[Problem 1.2]{aimlink}.

Motivated by the above question, in this paper we prove the following almost sure convergence for the number of the real roots inside $[-1,1]$, which could be viewed as a local strong law of large numbers.  As is well-known, the interval $[-1,1]$ plays an important role in the theory of Kac polynomials,  acting like a fundamental domain, since the distribution of the real roots of $p_n$ is essentially invariant under the symmetry $x\mapsto \frac 1 x$.  It is also classical that $\mathbb E N_n[-1,1]  = \frac {1}\pi \log n+o(\log n).$
 
 \begin{theorem} \label{t.main} Assume that $(\xi_j)_{j\ge 0}$ are iid  with zero mean,  unit variance, and  bounded $(2+\epsilon)^{th}$ moment for some $\epsilon>0$.   Then almost surely the following convergence holds:
\begin{equation*}
\lim_{n\to\infty} \frac{N_n([-1,1])}{\log n}=\frac 1\pi.
\end{equation*}
Furthermore, analogous results hold for $N_n[0,1]$ and $N_n[-1,0]$.
 \end{theorem}

Almost sure convergence  with SLLN flavors for real roots of random polynomials has recently attracted the attention of many authors.  In the Gaussian setting,  Ancona and Letendre \cite{AL}   established a SLLN for the   elliptic (also known as Kostlan) polynomials.  For random trigonometric functions,  a SLLN was established by Angst--Poly \cite{AP} when the iid random coefficients have a symmetric distribution and finite fourth moment.  For other models,  there are  recent partial results in this direction.  In  \cite{anw},  as consequences of their concentration estimates, Ager--Nguyen-Wang proved that $N_n(I)/\mathbb E[N_n(I)]\to 1$ for a large class of random functions (including elliptic polynomials and random Wey polynomials and random orthogonal polynomials) and $I$ is an appropriate subinterval of $\mathbb R$ that may capture the bulk of the real roots.   We note that in \cite{anw,AP} the results do not require  Gaussian distributions.

\subsection{Strategy}

Our method is inspired by ideas from harmonic analysis,  here to study  almost sure convergence of $N_n/\log n$ we will prove new maximal inequalities that enable us to reduce the proof to the simpler lacunary setting, where the desired convergence essentially follows from a recent concentration estimate  for $N_n$ by Can--Nguyen \cite{CN}.  As part of the proof of the maximal estimates, we will prove new small ball inequalities for $p_n$ that may also be of independent interest. The proof of the maximal inequalities will also involve a root repulsion argument inspired by prior joint work of the author with collaborators \cite{DHV,NNV}.  More details about our approach will be discussed in Section~\ref{s.outline}. See also Section~\ref{s.conclusion} for a related discussion of the limitation our methods when applying to the entire  $\mathbb R$.

\subsection{Notations and convention}

 $N_n(I)$ will denote the number of real roots of $p_n$ inside $I$ a subset of $\mathbb R$.  When it is clear from the context we might just drop the brackets from the notation, so  $N_n[0,1]$ will mean $N_n([0,1])$.

$[t]$ will denote the largest integer that does not exceed $t$.

We say that $A\ll B$ (equivalently $A=O(B)$) if there is an absolute constant $C>0$ such that $A\le CB$.  The dependence of $C$ on some parameters will be specified if needed.

We say that $A\approx B$ if $A\ll B$ and $B\ll A$.

\section{Outline of the proof of Theorem~\ref{t.main}}\label{s.outline}

Thanks to symmetries of the Kac polynomials (and of the proof), to show Theorem~\ref{t.main} it suffices to show that the following holds  almost surely 
\begin{eqnarray}\label{e.01}
\lim_{n\to\infty} \frac{N_n[0,1]}{\log n} &=& \frac 1{2\pi}.
\end{eqnarray}

Now,  to show \eqref{e.01},  we do not need to assume that $\xi_j$ are iid.  It suffices to assume that the distribution of $\xi_j$ \emph{does not depend on the degree $n$}\footnote{In other words,  $p_n(x)$ can be viewed as the $(n+1)$-th partial sum of a random power series $\sum_{j\ge 0}\xi_j x^j$.}, and these coefficients are independent with zero mean, unit variance, and uniformly bounded $(2+\epsilon)^{th}$ moments.   This will be assumed through the proof of \eqref{e.01}.

As a by-product of our estimates,   the inclusion/exclusion of the endpoints at $0$ and $1$ will not affect the left hand side of \eqref{e.01}. In fact, the following results hold: 
\begin{lemma}\label{l.near0-conv} Assume that the coefficients $\xi_j$ of $p_n$ are independent with  unit variance and uniformly bounded $(2+\epsilon)^{th}$ moments.
Then for any $C>1$ the following convergence holds almost surely,
\begin{eqnarray*}
\lim_{n\to\infty} \frac{N_n(-\frac 1 {C}, \frac 1 C)}{\log n} &=& 0.
\end{eqnarray*}
\end{lemma}
(Note that in Lemma~\ref{l.near0-conv} we do not need to assume that $\mathbb E \xi_j=0$.  The proof of Lemma~\ref{l.near0-conv} will be presented in Section~\ref{s.near0}.)

\begin{lemma}\label{l.near1-conv}
Assume that the coefficients $\xi_j$ of $p_n$ are independent with zero mean, unit variance, and uniformly bounded $(2+\epsilon)^{th}$ moments.
Then for any $C>1$ the following convergence holds almost surely,
\begin{eqnarray*}
\lim_{n\to\infty} \frac{N_n[1-\frac {C\log n}n,1]}{\log n} &=& 0.
\end{eqnarray*}
\end{lemma}
(The proof of Lemma~\ref{l.near0-conv} will be presented in Section~\ref{s.near1}.)

The simpler lacunary case of \eqref{e.01}  follows from a recent result of Can-Nguyen \cite{CN}.  
\begin{lemma} \label{l.lac} Assume that the coefficients $\xi_j$ of $p_n$ are independent with zero mean, unit variance, and uniformly bounded $(2+\epsilon)^{th}$ moments.
Let $1\le n_1<n_2<\dots$ be a   sequence of integers such that $\inf_k \frac{n_{k+1}}{n_k}>1$. Then the following holds almost surely
\begin{eqnarray*}
\lim_{k\to\infty} \frac {N_{n_k}[0,1]}{\log (n_k)} &=& \frac 1{2 \pi}.
\end{eqnarray*}
The same conclusion also holds for $[-1,0]$, $[1,\infty)$, $(-\infty,1]$, and $\mathbb R$. t
\end{lemma}

The rather straightforward derivation of Lemma~\ref{l.lac} from \cite{CN} will be presented later in Section~\ref{s.lac}.
The following maximal estimate  will be the key to reduce the proof of \eqref{e.01} to the lacunary setting:
\begin{lemma}\label{l.maximal} Assume that the coefficients $\xi_j$   are independent with zero mean, unit variance, and uniformly bounded $(2+\epsilon_0)^{th}$ moments.
Let $\epsilon>0$ and $c>1$ be deterministic constants independent of $n$.  Then there is a positive constant $c_1>0$ such that 
\begin{eqnarray*}
\Prob\Big(\max_{n\le k\le cn} \Big|N_k [0,1]-N_{n}[0,1]\Big| \ge   \epsilon\log n\Big) &\ll& n^{-c_1} .
\end{eqnarray*}
\end{lemma}
The proof of Lemma~\ref{l.maximal} will be presented in Section~\ref{s.maximal},  which uses ingredients in the proof of Lemma~\ref{l.near0-conv} and Lemma~\ref{l.near1-conv}.  

In the following,  we will show that Lemma~\ref{l.maximal} and Lemma~\ref{l.lac} imply Theorem~\ref{t.main}.  

Let $\epsilon>0$, and let $c=2$ and $n_j=2^j$.

For brevity, let $I:=[0,1]$.  Since $\sum_{j\ge 0} 2^{-c_1 j}<\infty$,  by  Borel-Cantelli's lemma and Lemma~\ref{l.maximal} the following holds: almost surely,  for $j$ large enough we have
\begin{eqnarray*}
\max_{n_j\le k\le n_{j+1}} |N_k (I)-N_{n_j}(I)|   &\le&   \epsilon\log {n_j},
\end{eqnarray*}
which implies
\begin{eqnarray*}
\max_{n_{j+1}\ge k\ge n_j} \frac {N_k(I)}{\log k}  &\le&   \epsilon +  \frac {N_{n_j}(I)}{\log n_j} .
\end{eqnarray*}
By taking  the supremum over $j$,  the following estimate holds almost surely for $j$ large enough,
\begin{eqnarray*}
\sup_{k\ge n_j} \frac {N_k(I)}{\log k}   &\le&     \epsilon + \sup(\frac {N_{n_j}(I)}{\log n_j}, \frac{N_{n_{j+1}}(I)}{\log n_{j+1}},\dots).
\end{eqnarray*}
By monotonicity and by  Lemma~\ref{l.lac},  almost surely we have
\begin{eqnarray*}
\limsup_{k\to\infty} \frac{N_k(I)}{\log k} &=& \lim_{j\to\infty} \sup_{k\ge n_j} \frac {N_k(I)}{\log k}\\
& \le& \epsilon+ \limsup_{j\to\infty} \frac{N_{n_j}(I)}{\log (n_j)}\\
&=&\epsilon+\frac 1 {2\pi}.
\end{eqnarray*}

By a similar argument (using Lemma~\ref{l.maximal})  the following holds almost surely,
\begin{eqnarray*}
\liminf_{k\to\infty} \frac{N_k(I)}{\log k}  & \ge&   -\epsilon + \liminf_{j\to\infty} \frac{N_{n_j}(I)}{\log (n_{j+1})},
\end{eqnarray*}
and using the fact that $\lim_{j\to\infty}\frac{\log n_{j+1}}{\log n_j} = 1$  and using Lemma~\ref{l.lac}, it follows that the following holds almost surely
\begin{eqnarray*}
\liminf_{k\to\infty} \frac{N_k(I)}{\log k} &\ge&  -\epsilon +\frac 1 {2\pi}.
\end{eqnarray*}
In the end, by sending $\epsilon\to 0$ we obtain the claim asserted in \eqref{e.01}.

\subsection{The lacunary case}\label{s.lac}

In this section, we will prove Lemma~\ref{l.lac}.  Namely, for every sequence $1\le n_1<n_2<\dots$ with $\inf_{k\ge 1}n_{k+1}/n_{k}=1+c>1$ the following holds almost surely:
\begin{eqnarray}\label{e.lac}
\lim_{k\to\infty} \frac{N_{n_k}[0,1]}{\log (n_k)} &=&  \frac 1 {2\pi}.
\end{eqnarray}
The proof is the same for other intervals, here we show the details for $[0,1]$. 

In Can--Nguyen \cite[Theorem 1.4, 1.5]{CN}, the authors show  that, for any $\epsilon>0$ there is some $c_\epsilon>0$ such that
\begin{eqnarray*}
\Prob(|N_n(I)-\mathbb E N_n (I)|\ge \epsilon \log n) &\ll_{\epsilon}&  e^{-c_\epsilon\sqrt{\log n}},
\end{eqnarray*}
where $I$ could be $\mathbb R$, $[0,1]$, $[1,\infty)$, $[-1,0]$, or $(-\infty,-1]$. In particular,
\begin{eqnarray*}
\Prob(|N_n(I)-\mathbb E N_n (I)|\ge \epsilon \log n) &\ll_{\epsilon}& (\log n)^{-2},
\end{eqnarray*}

Therefore, for any lacunary sequence $1\le n_1<n_2<\dots$ with $\inf \frac{n_{k+1}}{n_k}=1+c>1$, we have
\begin{eqnarray*}
\sum_{k\ge 1} \Prob(|N_{n_k}[0,1] - \mathbb E N_{n_k}[0,1]|>\epsilon \log {n_k}) &\ll_{c,\epsilon}& 1.
\end{eqnarray*}
Therefore, by Borel-Cantelli's lemma, almost surely the following holds  for $k$ large enough
\begin{eqnarray*}
|N_{n_k}[0,1]-\mathbb E N_{n_k}[0,1]| &\le& \epsilon \log {n_k}.
\end{eqnarray*}
Consequently,   the following holds almost surely
\begin{eqnarray*}
\limsup_{k\to\infty}\Big |\frac{N_{n_k}[0,1]}{\mathbb E N_{n_k}[0,1]} - 1 \Big | &\le& \limsup_{k\to\infty} \epsilon \frac{\log n_k}{\mathbb E N_{n_k}[0,1]} \ \ \ll\ \   \epsilon.
\end{eqnarray*}
Sending $\epsilon\to 0^+$ along a countable sequence, it follows that  \eqref{e.lac} holds almost surely.

\section{A small ball inequality for random polynomials} \label{s.smallball}

The goal of this section is to show the following small ball estimate,  which we will use in Section~\ref{s.near1}.  The estimate may be of independent interest. For convenience of notation, we let $V_n=V_n(x):=Var[p_n(x)]$.
\begin{theorem}\label{t.smallBE}Assume that $\xi_j$'s are independent with zero mean, unit variance, and uniformly bounded $(2+\epsilon)^{th}$ moments.   Let $c_0<1$ and $C_0>1$.  Then there is a positive constant $C_1>0$ such that the following holds for  $\frac 1 {C_0}\le x\le 1$ and any $\lambda$ satisfying
 $\lambda \gg \max(x^{c_0 n}V_n^{-1/2}, e^{-C_1V_n})$:
\begin{eqnarray}\label{e.smallBE}
\Prob(\frac 1{\sqrt{V_n}} |p_n(x)|\le \lambda) &\ll& \lambda.
\end{eqnarray}
\end{theorem}

We note that $V_n =\sum_{j=0}^n x^{2j} \approx (1-x+\frac 1n)^{-1}$ when $0\le x\le 1$, which will be useful in the application of Theorem~\ref{t.smallBE}.

It is clear that \eqref{e.smallBE} holds with no lower bound constraints for $\lambda$ if $\xi_j$ are furthermore Gaussian.   If $\xi_j$'s satisfy Cram\'er-type conditions or its extensions, Angst--Poly \cite{AP1} proved \eqref{e.smallBE} for $\lambda=n^{-s}$,  any $s>0$.  However, conditions of Cram\'er type are known to exclude the important cases when $\xi_j$'s  have discrete distributions, for example when $\xi_j=\pm 1$ with equal probability.\footnote{In this direction,  there are also   sublinear variants of Theorem~\ref{t.smallBE} (with  $\lambda^{1-c}$ on the right hand side but no normalization by $V_n$), see e.g. \cite[Lemma 7]{NNV} and also Appendix~\ref{app.a} of the current paper.} To avoid such Cram\'er conditions,  we will follow the approach in   \cite{DNN}  and \cite{KS} and prove a new estimate for the characteristic function of $p_n$ and then derive the conclusion of Theorem~\ref{t.smallBE} via the Fourier analytic argument of Hal\'asz \cite{Ha}.

In the rest of the section,  we assume that $\xi_j$ are independent with zero mean, unit variance, and uniformly bounded $(2+\epsilon)^{th}$ moments.  The implicit constants are  allowed to depend on $\epsilon$ and the uniform bound on the $(2+\epsilon)^{th}$ moments.

For convenience of future reference,  we write down a  corollary of Theorem~\ref{t.smallBE} that we will use in Section~\ref{s.near1}.
\begin{lemma}\label{l.smallBE}Let $C_0, C_1$ be positive constants.  Then there is a constant $c>0$ such that the following holds for all $\lambda \ge V_n^{-1/2}n^{-c}$ and all $x\in [1-\frac{C_0}{\log n},  1-\frac{C_1\log n}n]$:
\begin{eqnarray*}
\Prob(\frac 1{\sqrt{V_n}} |p_n(x)|\le \lambda) &\ll& \lambda.
\end{eqnarray*}
\end{lemma}
To see this consequence of Theorem~\ref{t.smallBE},  fix any $c_0\in (0,1)$,  and note that
\begin{eqnarray*}
V_n &\approx& (1-x+\frac 1n)^{-1} \ \ \gg\ \  \frac 1{C_0} \log n, \\
x^{c_0 n} &\le& (1-\frac{C_1\log n}n)^{c_0n}  \ \ \ll \ \  n^{-c_0'},
\end{eqnarray*}
for any $c_0' \le C_1c_0$.  By choosing $c_0'$ such that $C_1/(C_0c_0')$ is sufficiently large we will have $e^{-C_1 V_n}\ll n^{-c_0'}V_n^{-1/2}$.  In the end,  we simply choose $0<c<c_0'$ to obtain Lemma~\ref{l.smallBE} from Theorem~\ref{t.smallBE}.

\subsection{Proof of Theorem~\ref{t.smallBE}, step 1: a characteristic function bound}
We first prove a bound for the characteristic function $\Phi$ of $\frac1{\sqrt V_n}p_n(x)$. Let $i=\sqrt{-1}$, then define
\begin{eqnarray*}
\Phi(w) &:=& \mathbb E[\exp(2\pi i w V_n^{-1/2}p_n(x))].
\end{eqnarray*}

\begin{lemma} \label{l.charb-gen}
Let $c_0<1$ and $C_0>1$. Then there are positive constants $C_1, C_2$ such that for any $\frac 1{C_0}\le x\le 1$ and any   $w\in \mathbb R$ satisfying
\begin{eqnarray*}
|w| &\le& \frac 1 {C_2} (1-x+\frac 1n)^{-1/2} x^{-c_0 n}
\end{eqnarray*}
we have
\begin{eqnarray*}
|\Phi(w)| &\le& \exp(-C_1\min(V_n, w^2)).
\end{eqnarray*}
\end{lemma}

\proof[Proof of Lemma~\ref{l.charb-gen}] For convenience of notation, let $X_j=\frac{\xi_j x^j}{\sqrt V_n}$, then $\frac 1{\sqrt{V_n}}p_n(x) =\sum_{j=0}^n X_j$.  

Without loss of generality, assume $c_0\ge 0$. 

Given a random variable $\xi$, consider the $\xi$-norm (see \cite{TV08} for more background)
\begin{eqnarray*}
\|w\|_\xi &:=& (\mathbb E \|w(\eta_1-\eta_2)\|_{\mathbb R/\mathbb Z}^2)^{1/2},
\end{eqnarray*}
where $\eta_1$ and $\eta_2$ are two iid copies of $\xi$, and $\|.\|_{\mathbb R/\mathbb Z}$ denotes the distance to the nearest integer.
Then we have the elementary bound (see e.g. \cite{TV08})
\begin{eqnarray*}
|\mathbb E[e^{2\pi i\xi x}]| &\le& \frac 1 2 + \frac 1 2|\mathbb E[e^{2\pi i\xi x}]|^2 \\
&=& \mathbb E [\frac 1 2 + \frac 1 2 \cos(2\pi(\eta_1-\eta_2)x)]\\
&\le& e^{-c\|x\|_{\xi}^2}
\end{eqnarray*}
where $c>0$ is an absolute constant.  Therefore
\begin{eqnarray*}
|\Phi(w)| &=& \prod_{j=0}^n |\mathbb E[\exp(2\pi i w X_j)]| \\
&\le& \exp(-c\sum_{j=0}^n \|w\|_{X_j}^2).
\end{eqnarray*}
Thus, it suffices to prove the following lower bound
\begin{eqnarray*}
\sum_{j=0}^n \|w\|_{X_j}^2 &\gg& \min(w^2, V_n).
\end{eqnarray*}
Now, using the fact that $\xi_j$   have zero mean,  unit variance,  and uniformly bounded $(2+\epsilon)^{th}$ moments,  we can find a positive constant  $C_0$ that may depend on $\epsilon$ (and the uniform moment bound) such that $P(\frac 1 C_0 \le |\eta_{j,1}-\eta_{j,2}|\le C_0) \ge \frac 1{C_0}$ for all $j=0,1,\dots,n$,  here $\eta_{j,1}$ and $\eta_{j,2}$ are iid copies of $\xi_j$.\footnote{For a proof of this estimate,  see  \cite[Lemma 9.1]{TV}.}  It follows that
\begin{eqnarray*}
\sum_{j=0}^n \|w\|_{X_j}^2 
&=&  \sum_{j=0}^n \mathbb E \Big[\|\frac{wx^j}{\sqrt{V_n}}(\eta_{j,1}-\eta_{j,2})\|_{\mathbb R/\mathbb Z}^2\Big]  \\
&\gg&   \sum_{j=0}^n \inf_{\frac 1 C_0 \le |\eta|\le C_0} \|\frac{wx^j\eta}{\sqrt{V_n}}\|_{\mathbb R/\mathbb Z}^2 .
\end{eqnarray*}

We now consider those $j$'s such that 
\begin{eqnarray*}
|w|x^j  &<& \frac {\sqrt{V_n}}{2C_0}.
\end{eqnarray*}
Note that this property will ensure that $|\frac{wx^j\eta}{\sqrt{V_n}}|\le 1/2$  (for those $\eta$ inside the infimum range) and therefore $\|\frac{wx^j\eta}{\sqrt{V_n}}\|_{\mathbb R/\mathbb Z}=\frac{|wx^j\eta|}{\sqrt{V_n}}\gg \frac{|w|x^j}{\sqrt{V_n}}$.

We divide the consideration into two cases.

If  $|w| < \frac {\sqrt{V_n}}{2C_0}$ then all $j\ge 0$ are good,  and we will have
\begin{eqnarray*}
\sum_{j=0}^n \|w\|_{X_j}^2 
&\gg&  \sum_{j=0}^n (\frac 1{\sqrt{V_n}} |w|x^j)^2  \\
&\gg& w^2V_n^{-1}\sum_{j=0}^n x^{2j}\\
&\gg& w^2.
\end{eqnarray*}

If $|w| \ge \frac {\sqrt{V_n}}{2C_0}$ then the smallest $j$ we can take is $j_0$, 
defined by 
\begin{eqnarray*}
j_0 &:=&  1+[\frac{\log(\frac{2C_0|w|}{\sqrt{V_n}})}{\log(1/x)}].
\end{eqnarray*}
Note that the given hypothesis on $w$ implies
\begin{eqnarray*}
|w|x^{c_0 n} &\le& \frac 1 {C_2} (1-x+\frac 1 n)^{-1/2}\\
&\ll& \frac {\sqrt{V_n}}{C_2},
\end{eqnarray*}
therefore by choosing $C_2>0$ large enough we can see that $j_0\le c_0 n$.

Since $x\ge \frac 1{C_0}>0$,  we also have $x^{j_0} \gg x^{j_0-1} \gg \frac {\sqrt {V_n}}{|w|}$.  Consequently,
\begin{eqnarray*}
\sum_{j=0}^n \|w\|_{X_j}^2 
&\gg&  \sum_{j=j_0}^n (\frac 1{\sqrt{V_n}} |w|x^j)^2  \\
&\gg& \frac{1}{V_n}x^{2j_0}|w|^2\sum_{k=0}^{n-c_0n}x^{2k}\\
&\gg& V_n.
\end{eqnarray*}
This completes the   proof of Lemma~\ref{l.charb-gen}.
\endproof

\subsection{Proof of Theorem~\ref{t.smallBE}, step 2: proof of the small ball estimates}
We now discuss the rest of the proof of Theorem~\ref{t.smallBE}. 

Without loss of generality assume $c_0\ge 0$. We will use the constants $C_1,C_2$ that appear  in the conclusion of Lemma~\ref{l.charb-gen}.

Recall that $\Phi$ denotes the characteristic function of $V_n^{-1/2}p_n(x)$.  

Choose $\varphi\in \mathbb C^\infty(\mathbb R)$ such that $\varphi$ is nonnegative and $\varphi\ge 1$ in $(-1,1)$ and the support of its Fourier transform $\widehat \varphi$ is inside $(-C,C)$ for some   absolute constant $C>0$.   As usual, we define the Fourier transform
\begin{eqnarray*}
\widehat \varphi(\xi) &=& \int_{\mathbb R} e^{-2\pi i x\xi}\varphi(x)dx.
\end{eqnarray*}
For instance, consider $\rho$ a function whose Fourier transform is a nonzero bump function supported inside $(-1,1)$, then $|\rho|^2$ is continuous, nonnegative with frequency support in $(-2,2)$ and one may further dilate and multiply $|\rho|^2$ to get $\varphi$ with the desired lower bound $\varphi\ge 1$ on $(-1,1)$.  (The dilation may enlarge the compact frequency support.)

By monotonicity of the distribution function and by replacing $\lambda$ by $C'\lambda$ for a large   constant $C'>0$ if necessary, we may assume without loss of generality that $\lambda \ge CC_2 \max(x^{c_0 n}(1-x+\frac 1n)^{1/2}, e^{-C_1V_n})$,  so in particular Lemma~\ref{l.charb-gen} is applicable when $|w|\le \frac C\lambda$.

We obtain
\begin{eqnarray*}
\Prob(\frac 1{\sqrt{V_n}} |p_n(x)|\le \lambda) &\le& \mathbb E [\varphi(\frac 1{\lambda\sqrt{V_n}} p_n(x))]\\
&=&  \mathbb E [\int_{\mathbb R} \exp(2\pi i w \frac {p_n(x)}{\lambda \sqrt{V_n}}) \widehat \varphi(w)dw]\\
&=&  \int_{\mathbb R} \Phi(\frac{w}\lambda) \widehat \varphi(w)dw\\
&=& \lambda \int_{|w|\le C/\lambda} \Phi(w) \widehat \varphi(\lambda w)dw .
\end{eqnarray*}
Using Lemma~\ref{l.charb-gen}, we obtain
\begin{eqnarray*}
\Prob(\frac 1{\sqrt{V_n}} |p_n(x)|\le \lambda) 
&\ll&
 \lambda \int_{|w|\le C/\lambda} e^{-C_1\min(V_n,w^2)}dw
\end{eqnarray*}

We will split the integral into two regions $I=\{|w|\le \sqrt{V_n}\}$ and $II=\{\sqrt{V_n}<|w|\le C/\lambda\}$ and estimate them separately.  Clearly,
\[
|\lambda \int_I e^{-C_1 w^2}dw | \ \  \le  \ \ 2|\lambda \int_0^\infty e^{-C_1 w^2}dw| \ \ \ll \ \   \lambda,
\]
\[
|\lambda \int_{II} e^{-C_1 V_n}dw | \ \ \ll \ \  (\lambda)(\frac 1 \lambda)  e^{-C_1 V_n} \ \  \ll \ \  \lambda.
\]
This completes the proof of Theorem~\ref{t.smallBE}. \endproof

\section{Maximal estimates and   convergence for small roots}\label{s.near0}

In this section,  we consider maximal estimates involving the number of real roots inside $(-\frac 1 C, \frac 1 C)$, where $C>1$ is an  absolute constant.  As a corollary of our estimates,  we will  obtain the following result.
\begin{lemma}\label{l.near0-conv-alt} Assume that the coefficients $\xi_j$ of $p_n$ are independent with $\mathbb E|\xi_j|=O(1)$ and such that
for some $0<c_0<\infty$ and $0<q_0<1$ we have
\begin{eqnarray}\label{e.smallball}
\Prob(|\xi_j|<c_0) &\le& q_0.
\end{eqnarray}
Then for any $C>1$ the following convergence holds almost surely,
\begin{eqnarray*}
\lim_{n\to\infty} \frac{N_n(-\frac 1 {C}, \frac 1 C)}{\log n} &=& 0.
\end{eqnarray*}
\end{lemma}
 
We first discuss how Lemma~\ref{l.near0-conv-alt} implies Lemma~\ref{l.near0-conv},  which we mentioned in Section~\ref{s.outline}.  For this implication, it suffices to show that the condition \eqref{e.smallball} would be satisfied if we assume that the random coefficients $\xi_j$ have unit variance and uniformly bounded $(2+\epsilon)^{th}$ moments.  Indeed,  let $\widetilde \xi_j=\widetilde \xi_j-\mathbb E \widetilde \xi_j$,  then the random  variables $\widetilde \xi_j$  have zero mean,  unit variance and uniformly bounded $(2+\epsilon)^{th}$ moments.   It follows that there exists $C>1$ (uniform over $j$) such that 
\begin{eqnarray*}
\Prob(C^{-1} \le |\widetilde \xi_j- \widetilde \xi_j'|\le C) &\ge& C^{-1},
\end{eqnarray*}
where $\widetilde \xi'_j$ is constructed in the same way from an independent copy $\xi'_j$ of $\xi_j$.  Now it is clear that $\widetilde \xi_j -\widetilde \xi_j' = \xi_j-\xi'_j$.  Therefore,
\begin{eqnarray*}
\Big(\Prob(|\xi_j| < \frac 1 {3 C})\Big)^2 &=& \Prob(|\xi_j|< \frac 1{3 C}, |\xi_j'| < \frac 1 {3 C})  \\
& \le& \Prob(|\xi_j-\xi_j'|<C^{-1}) \\
&\le& 1- C^{-1},
\end{eqnarray*}
so we may take $c_0=\frac 1 {3C}$ and $q_0=(1-\frac 1 C)^{1/2}$.

The key ingredient in our proof of Lemma~\ref{l.near0-conv-alt} is the following maximal estimate:

\begin{lemma} \label{l.near0-max}
Assume that the coefficients $\xi_j$ of $p_n$ are independent with $\max_j \mathbb E|\xi_j| =O(1)$, and  for some $0<c_0<\infty$ and $0<q_0<1$ we have
\begin{eqnarray*}
\Prob(|\xi_j|<c_0) &\le& q_0.
\end{eqnarray*}

Then, given any finite constant $C_1>1$, there exists a positive constant $c_2>0$ such that the following estimate holds uniformly over $n$ and $t$:
\begin{eqnarray}\label{e.near0-max}
\Prob (\max_{0\le j\le n} N_j(-\frac 1 {C_1},\frac 1 {C_1})\ge t) &\ll& e^{-c_2t}.  
\end{eqnarray}
Furthermore,  we may take $c_2$ comparable to $\log(1/q_0)$ if $C_1$ is sufficiently large (depending on $q_0$). 

\end{lemma}

\begin{remark}
Generally speaking, exponential decaying estimates are the best we can achieve in \eqref{e.near0-max}, and the threshold value $c_2=\log(1/q_0)$ is sharp regardless of how large we take $C_1$.  To see this, suppose that $\xi_j$ are iid with  an atom at $0$, i.e.  $q_1:=\Prob(|\xi_j|=0)$ is strictly between $0$ and $1$, then for any $t\ge 0$
\begin{eqnarray*}
\Prob(N_n(\{0\}) \ge t) &\ge& \Prob(\xi_0= \dots=\xi_{[t]-1}=0) \\
&\gg& q_1^{[t]} \\
&\gg& e^{-t\log(1/q_1)}.
\end{eqnarray*}
Furthermore, if  $|\xi_j|$ has a discrete distribution with support on $\{-c,0,c\}$ for some $c>0$ then it is clear that $q_0=q_1$,  thus $c_2=\log(1/q_0)$ is the best exponent regardless of how large  $C_1$ is chosen.
\end{remark}

\subsection{From maximal estimates to almost sure convergence: Proof of Lemma~\ref{l.near0-conv-alt} using Lemma~\ref{l.near0-max}}
We will use the following basic lemma (which we will use again in Section~\ref{s.near1}).
\begin{lemma}\label{l.decay->pwconv}  
Let $(X_n)_{n\ge 1}$ be a sequence of random variables and assume that for any $\epsilon>0$ there is $C_0>0$ and $C_1>1$ such that the following estimate holds uniformly over   $n$,
\begin{eqnarray*}\label{e.maximal}
\Prob(\max_{n\le j\le C_1 n}|X_j| \ge \epsilon) &\ll_\epsilon& (\log n)^{-(2+C_0)}.
\end{eqnarray*}
Then almost surely it holds that
\begin{eqnarray*}\label{e.pwconv}
\lim_{n\to\infty} X_n &=& 0.
\end{eqnarray*}
\end{lemma} 
\proof[Proof of Lemma~\ref{l.decay->pwconv}]
Let $\epsilon>0$ be given. By an union bound, for any $m\ge 1$ integer we have
\begin{eqnarray*}
\Prob(\max_{n\le j\le C_1^m n}|X_j| \ge \epsilon) &\ll_{\epsilon,m}& (\log n)^{-(2+C_0)},
\end{eqnarray*}
therefore without loss of generality we may assume $C_1=2$.

Now, the following holds uniformly over all $k$ sufficiently large,
\begin{eqnarray*}
\Prob(\max_{2^k\le j\le 2^{k+1}} |X_j|\ge \epsilon) &\ll& k^{-(2+C_0)}.
\end{eqnarray*}
Thus, by an union bound,
\begin{eqnarray*}
\Prob(\sup_{j\ge 2^k} |X_j| \ge \epsilon) &\ll&  \sum_{m\ge k} m^{-(2+C_0)} \\
&\ll& k^{-(1+C_0)}.
\end{eqnarray*}
Thus, using Borel-Cantelli's lemma, it follows that 
\begin{eqnarray*}
\lim_{k\to\infty} \sup_{j\ge 2^k} |X_n| &\le&  \epsilon,
\end{eqnarray*}
almost surely. Note that the sequence $\sup_{j\ge n} |X_j|$ is non-increasing with respect to $n$, therefore the limit on the left hand side exists.  In particular, we have
\begin{eqnarray*}
\limsup_{n\to\infty} |X_n|  &=& \lim_{n\to \infty} \sup_{j\ge n} |X_j| \\
&\equiv& \lim_{k\to\infty} \sup_{j\ge 2^{k}}  |X_j| \\
&\le& \epsilon. 
\end{eqnarray*}
Since this holds for any $\epsilon>0$ we obtain
\begin{eqnarray*}
\lim_{n\to\infty} |X_n| &=& 0,
\end{eqnarray*}
which completes the current proof.
\endproof
We now discuss how Lemma~\ref{l.near0-conv-alt} follows from Lemma~\ref{l.near0-max}, with the aid of Lemma~\ref{l.decay->pwconv}. Let $C>1$. By Lemma~\ref{l.near0-max}, we have
\begin{eqnarray*}
\Prob(\max_{n\le j\le 2n} N_j(-\frac 1 C, \frac 1 C) \ge \epsilon \log n) &=& O(n^{-C'\epsilon})
\end{eqnarray*}
for some $C'>0$. Since $  \log(2n) >\log n>0$ (for $n\ge 2$), we obtain
\begin{eqnarray*}
\Prob(\max_{n\le j\le 2n} \frac{N_j(-\frac 1 C, \frac 1 C)}{\log j} \ge \epsilon) &=& O(n^{-C'\epsilon}),
\end{eqnarray*}
and so by Lemma~\ref{l.decay->pwconv} we obtain the desired convergence in Lemma~\ref{l.near0-conv-alt}.


\subsection{Proof of Lemma~\ref{l.near0-max}}

Without loss of generality, we may assume that $0< t\le n$ and $0<c_0<1$.  Define
\begin{eqnarray*}
B_0 &=& \{|\xi_0|\ge c_0\},\\
B_k &=&\{|\xi_0|<c_0,\dots,|\xi_{k-1}|<c_0, |\xi_k|\ge c_0\}, \ \ \dots
\end{eqnarray*}
for $0\le k \le n$.  Then
\begin{eqnarray*}
\Prob(B_k) &\le& q_0^k, 
\end{eqnarray*}
and
\begin{eqnarray*}
\Prob((B_0\cup\dots\cup B_n)^c) &=& \Prob(\max_{0\le j\le n}|\xi_j|<c_0) \\
&\le& q_0^{n+1} \\
&=& O(e^{-t\log(1/q_0)}).
\end{eqnarray*}
Thanks to the last estimate,  it follows that for the current proof we may restrict our attention to   $B_0\cup \dots \cup B_n$.

Let $r=\frac 1{C_1}$ and $R\in (r,1)$,  to be chosen later that may depend on $C_1$ and $q_0$. 

We will show that, for any $0\le k,j\le n$,  on the event $B_k$ the following holds
\begin{eqnarray*}
N_{j}(-\frac 1 {C_1}, \frac 1{C_1})  
&\le& k+\frac {\log(M_k/k!)}{\log (R/r)} + \frac{\log(1/c_0)}{\log (R/r)},\\
M_k &:=&\sum_{m=k}^n |\xi_m| R^{m-k} \frac{m!}{(m-k)!}.
\end{eqnarray*}

Note that on $B_k$ we have $M_k\ge k!|\xi_k| \ge k!c_0$,  thus the right hand side of the above estimate is at least $k$.  Therefore,  if $j \le k$ the above estimate follows from the trivial estimate $N_j\le j$.

For $k< j\le n$,  on $B_k$ we use the mean value theorem and Jensen's formula to estimate
\begin{eqnarray*}
N_{j}(-\frac 1 {C_1}, \frac 1{C_1}) &\le& k+\frac{\log \frac{\sup_{|z|=R}|p^{(k)}_j(z)|}{|p_j^{(k)}(0)|}}{\log \frac R r} \\
&\le& k+\frac {\log(M_k/k!)}{\log (R/r)} + \frac{\log(1/c_0)}{\log (R/r)},
\end{eqnarray*}
again giving us the above estimate for $N_j(-\frac 1{C_1},\frac 1{C_1})$.

It follows that, on $B_0\cup\dots\cup B_n$ we have
\begin{eqnarray}\label{e.maxJensen}
&& \max_{j\le n} N_j(-\frac 1 {C_1},\frac 1 {C_1}) \\
\nonumber &\le& \frac{\log(1/c_0)}{\log (R/r)}+
\sum_{0\le k\le n} k1_{B_k} + \frac 1{\log(R/r)}\sum_{0\le k\le n} \log(M_k/k!)1_{B_k}.
\end{eqnarray} 
As our choice of $R$ will satisfy the lower bound $\log(R/r) \gg 1$, the first term in \eqref{e.maxJensen} is $O(1)$ therefore  we may safely ignore it in the current proof.

We now estimate the impact of the second term in the estimate \eqref{e.maxJensen} towards the left hand side of \eqref{e.near0-max}. We have
\begin{eqnarray*}
\Prob(\sum_k k1_{B_k} \ge t) &=& \sum_{k\ge t}\Prob(B_k)\\
&=& O(q_0^{-t}) \\
&=& O(e^{-t \log (1/q_0)}).
\end{eqnarray*}
For the contribution of the third (and last) term  in the estimate \eqref{e.maxJensen}, we note that
\begin{eqnarray*}
\mathbb E [M_k] &\ll& \sum_{j=k}^\infty \frac{j!}{(j-k)!}R^{j-k}\ \ = \ \   k!(1-R)^{-(k+1)}.
\end{eqnarray*}
Therefore, 
\begin{eqnarray*}
\Prob(\frac{M_k}{k!} \ge \alpha)  &\ll& \alpha^{-1}(1-R)^{-(k+1)}.
\end{eqnarray*}
Consequently,
\begin{eqnarray*}
&& \Prob(\frac 1{\log (R/r)} \sum_{k=0}^n \log(M_k/k!)  1_{B_k}\ge t) \\
&\le& \sum_{k=0}^n \min\Big(\Prob(M_k/k! \ge  e^{t\log(R/r)}),\Prob(B_k)\Big)\\
&\ll&  \sum_{k=0}^n \min\Big(e^{-t\log(R/r)}(1-R)^{-(k+1)},q_0^{k}\Big).
\end{eqnarray*}
The last estimate is a two-sided geometric sum that can be estimated by standard arguments, as follows.
Let $\kappa_0\in \mathbb R$ be such that $(q_0(1-R))^{\kappa_0} = e^{-t\log(R/r)}$, in other words
\[
\kappa_0=\frac{-t\log(R/r)}{\log(q_0(1-R))}\ge 0.
\]
It follows that (below the implicit constants may depend on $R,  r, q_0$, but not on $t$)
\begin{eqnarray*}
&&\Prob(\frac 1{\log (R/r)} \sum_{k=0}^n \log(M_k/k!)  1_{B_k}\ge t)\\
&\le& \sum_{k\le \kappa_0} e^{-t\log(R/r)}(1-R)^{-(1+k)}+\sum_{k>\kappa_0}q_0^k  \\
&\ll& e^{-t\log(R/r)}(1-R)^{-(\kappa_0+1)} + q_0^{\kappa_0} \\
&\ll& q_0^{k_0} \\
&=& e^{-tC_3}, \  \ \ \ C_3:=\frac{\log(1/q_0)\log(R/r)}{\log(1/q_0)+\log(1/(1-R))}>0.
\end{eqnarray*}
There may be an optimal choice for $R\in (\frac 1 {C_1},1)$, but for our purpose we may take for instance $R=\frac 1 2(1+\frac {1}{C_1})\in (\frac 1 {C_1},1)$. Collecting estimates, we obtain
\begin{eqnarray*}
\Prob (\max_{0\le j\le n} N_j(-\frac 1 {C_1},\frac 1 {C_1})\ge t) &=& O(e^{-C_3t}) + O(e^{-t\log(1/q_0)}).  
\end{eqnarray*}
 The positive constant $c_2$ in the desired estimate can be taken to be the minimum of $C_3$ and $\log(1/q_0)$.  
This completes the proof of Lemma~\ref{l.near0-max}.

\section{Maximal estimates and  convergence near $1$}\label{s.near1}

In this section,  we will assume that the coefficients $\xi_j$ are independent with zero mean,  unit variance,  and uniformly bounded $(2+\epsilon_0)^{th}$ moments.

The goal of this section is to prove the following maximal inequality for the number real zeros inside the interval $J:=[1-\frac {C\log n}n,1]$,  where $C$ is an absolute constant.  

\begin{lemma} \label{l.near1-max} Let $\epsilon>0$ and $c>1$ and $C>0$. Then there is a constant $c'>0$ such that  
\begin{eqnarray*}
\Prob(\max_{n\le m\le cn} N_m[1-\frac{C\log n}n, 1]>\epsilon \log n) &\ll& n^{-c'}.
\end{eqnarray*}
\end{lemma}

\begin{remark} Since $\log(cn)\approx \log n$,  it follows from Lemma~\ref{l.near1-max}
that
\begin{eqnarray*}
\Prob(\max_{n\le m\le cn} \frac{N_m[1-\frac{C\log m}m, 1]}{\log m}>\epsilon) &\ll& n^{-c'}.
\end{eqnarray*}
It is clear that  this estimate and  Lemma~\ref{l.decay->pwconv} implies Lemma~\ref{l.near1-conv}.
\end{remark}

\subsection{An auxiliary estimate}
We first prove an estimate for iterated integrals that we will use in the proof of Lemma~\ref{l.near1-max}.  Similar estimates were  considered in Maslova's work \cite{M,M1},  and reappeared in recent studies \cite{NV1, CN}.  

\begin{lemma}\label{l.iterated-bound} For some absolute constant $C>0$ the following holds: Suppose that $1\le k \le n$ and $0\le x\le y\le 1$, then
\begin{eqnarray}\label{e.iterated-bound}
&&\int_{x}^{y}\int_{x}^{y_1}\dots \int_{x}^{y_{k-1}} \mathbb E |p_n^{(k)}(y_k)|^2 dy_k\dots dy_1 \\
\nonumber &\le& C\min\Big(y^k (\frac n2)^{k+1},     k! \sqrt k   (\frac{1}{2(1 -y)})^{k+1}  \Big).
\end{eqnarray}
Similarly, we also have
\begin{eqnarray}\label{e.iterated-bound-dual}
&&\int_{x}^{y}\int_{x_1}^{y}\dots \int_{x_{k-1}}^{y} \mathbb E |p_n^{(k)}(x_k)|^2 dx_k\dots dx_1 \\
\nonumber &\le& Cn^{2k+1}\frac{(1-x)^k}{k!}.
\end{eqnarray}

\end{lemma}
\begin{remark} It is important that the constant $C$ does not depend on $k$,  since $k$ will depend on $n$ in the applications of this estimate.
\end{remark}
\proof
For brevity, within this proof let $LHS$ denote the left hand side of \eqref{e.iterated-bound}. Without loss of generality we may even let $x=0$ in \eqref{e.iterated-bound}.
Since the coefficients $\xi_j$ are independent with zero mean and unit variance,  we have
\begin{eqnarray*}
LHS &\le &  \sum_{j=0}^{n-k} (\frac {(j+k)!}{j!})^2 \int_0^y \dots \int_0^{y_{k-1}} y_k^{2j}  dy_k\dots dy_1 \\
&=&  \sum_{j=0}^{n-k}  (\frac {(j+k)!}{j!})^2 \frac{(2j)!}{(2j+k)!} y^{2j+k}.
\end{eqnarray*}

Using Stirling's formula $n!\sim \sqrt{2\pi n} (n/e)^n$, we may estimate,  
\begin{eqnarray*}
(\frac {(j+k)!}{j!})^2 \frac{(2j)!}{(2j+k)!}  &\ll& (1+j)^{-\frac12} e^{-k}\frac{2^{2j}(j+k)^{2j+2k+1}}{(2j+k)^{2j+k+\frac 12}} \\
&=&(1+j)^{-\frac12} 2^{-k} \frac{(j+k)^{k+1}}{(2j+k)^{1/2}} \Big(e^{-k}(1+\frac k{2j+k})^{2j+k}\Big)\\
&\le& (1+j)^{-\frac12} 2^{-k} (j+k)^{k+\frac 1 2}.
\end{eqnarray*}
In particular, we obtain
\begin{eqnarray*}
LHS &\ll&  y^k  \sum_{0\le j\le n-k}  (j+1)^{-1/2}  (\frac n2)^{k+\frac 1 2} \\
&\ll&  y^k (\frac n2)^{k+1} ,
\end{eqnarray*}
which is the first part of the asserted estimate. 

For the second estimate, we proceed by
\begin{eqnarray*}
LHS  &\ll&\sum_{0\le j\le n-k}  (1+j)^{-\frac12} 2^{-k} (j+k)^{k+\frac 1 2}y^{2j+k} \\
&\le&    2^{-k}  \sum_{j\ge 0}\sqrt{\frac{j+k}{j+1}} \frac{(2j+2k)!}{(2j+k)!}y^{2j+k}\\
&\le& 2^{-k}\sqrt k \sum_{m\ge 0} \frac{(m+k)!}{m!}y^{m} \\
&=& 2^{-k}\sqrt k k! (1 -y)^{-(k+1)}  .
\end{eqnarray*}

We now turn to the proof of \eqref{e.iterated-bound-dual}.
For brevity, within this proof let $LHS$ denote the left hand side of \eqref{e.iterated-bound-dual}. Without loss of generality we may even let $y=1$ in \eqref{e.iterated-bound-dual}.
Since the coefficients $\xi_j$ are independent with zero mean and unit variance,  we have
\begin{eqnarray*}
LHS &\le &  \sum_{j=0}^{n-k} (\frac {(j+k)!}{j!})^2 \int_x^1 \dots \int_{x_{k-1}}^1 1 dx_k\dots dx_1 \\
&=&  \sum_{j=0}^{n-k}  (\frac {(j+k)!}{j!})^2 \frac{(1-x)^k}{k!}\\
&\le&   n^{2k+1}\frac{(1-x)^k}{k!}.
\end{eqnarray*}
\endproof

\subsection{Proof of Lemma~\ref{l.near1-max}}

Our strategy is to divide $[1-\frac {C\log n}n,1]$ into smaller intervals and prove the analogous estimate for each interval and then use the triangle inequality to complete the proof.  Among the intervals,  the end interval shall be of the form $[1-\frac {\log n}{C'n},1]$ where $C' \in(\frac 1 C,\infty)$ is a large constant (depending on $c$ and $\epsilon$, in particular we will take $C'>e/(2\epsilon)$). The remaining interval $[1-\frac {C\log n}n, 1-\frac{\log n}{C'n}]$ will be further divided into  subintervals and the number of subintervals will be $O(1)$ (independent of $\epsilon$).  The endpoints for each such subinterval  will be close in the pseudo-hyperbolic distance $d(x,y):=\frac{|x-y|}{|1-xy|}$.  More specifically,   we'll divide it into $L$ subintervals $[\alpha_j,\alpha_{j+1}]$  where 
\begin{eqnarray}\label{e.alphaj}
\alpha_j &:=& 1-C(C'C)^{-j/L} \frac{\log n}n
\end{eqnarray}
and $j=0,1,\dots,L$.   We will choose $L$ large enough depending on $C,C'$. It is clear that the hyperbolic distance between the endpoints  becomes (uniformly) small if $L$ is large:
\begin{eqnarray*}
d(\alpha_j,\alpha_{j+1}) &\le& \frac{|\alpha_{j+1}-\alpha_j|}{|1-\alpha_{j+1}|} \\
&\ll_{C,C'}& \frac{1}{L}.
\end{eqnarray*}
We stress that the implicit constants in the above display shall depend only on $C,C'$ (and not $n$ or $L$).  In particular,  they are also uniform over $j$.

For each sub-interval $[\alpha_j,\alpha_{j+1}]$, we will use the following estimate.

\begin{lemma} \label{l.near1-max-sub} Let $c>1$ and let $C_0,C_0'$ be positive constants.  Then for any $c_0>0$ sufficiently small the following holds.  Let $L>0$ that may depend on $n$ and assume   $0\le x\le y\le 1$ such that $|y-x|\le \frac 1{L} |1-y|$ and $\frac{C_0\log n}{n}\le |1-y| \le \frac {C_1}{\log n}$. Then  for any $1\le k\le n$  
\begin{eqnarray*}
&&\Prob(\max_{n\le m\le cn} N_m[x,y] \ge k)\\
& \ll& n^{-c_0} + n^{2c_0} \sqrt{k}(2L)^{-k}+n^{2c_0}y^{2n}.
\end{eqnarray*}
Here, the implicit constant may depend on $c_0,c$ and the implicit constants in the assumptions. 
\end{lemma}

For the end interval (with right endpoint at $1$) we will use the following estimate.

\begin{lemma} \label{l.near1-max-end} Let $c>1$ and let $C_0,C_0'$ be positive constants..  Then  for any $c_0>0$ sufficiently small the following holds.  Assume   $0\le x\le 1$ such that $\frac{C_0\log n}n \le |1-x|\le \frac {C_1}{\log n}$. Then  for any $1\le k\le n$ 
\begin{eqnarray*}
&&\Prob(\max_{n\le m\le cn} N_m[x,1] \ge k) \\
&\ll&  n^{-c_0} +   n^{2c_0}  (\frac{cen(1-x)}{k})^{2k+1}  +n^{2c_0} x^{2n}.
\end{eqnarray*}
Here, the implicit constant $C$ may depend on $c_0,c$ and the implicit constants in the assumptions. 
\end{lemma}

We first discuss the treatment of $[1-\frac {C\log n}n, 1-\frac{\log n}{C'n}]$ using Lemma~\ref{l.near1-max-sub}.  Let $1\le L \le \epsilon\log n$ and perform the division of this interval into $[\alpha_j,\alpha_{j+1}]$ with $\alpha_j$ defined by \eqref{e.alphaj}.  Let
\begin{eqnarray*}
k&=&1+[\frac 1L\epsilon \log n].
\end{eqnarray*}
Then as discussed above there is some implicit constant $C_1>0$ that may depend only on $C,C'$ such that
\begin{eqnarray*}
\frac{\alpha_{j+1}-\alpha_j}{1-\alpha_{j+1}} &\le& \frac{C_1}{L}.
\end{eqnarray*}
Note that $\frac {\log n}{C'n} \le 1-\alpha_{j+1} \le \frac {C\log n}{n}\ll \frac{1}{\log n}$ uniformly over $j$.  Thus,  applying   Lemma~\ref{l.near1-max-sub} we obtain, for $c_0>0$ sufficiently small (uniform over $j$),
\begin{eqnarray*}
&& \Prob(\max_{n\le m \le cn} N_m[1-\frac {C\log n}n, 1-\frac{\log n}{C'n}]> \epsilon \log n)  \\
&\le& \sum_{j=0}^L \Prob(\max_{n\le m \le cn} N_m[\alpha_j,  \alpha_{j+1}]\ge k) \\
&\ll& Ln^{-c_0} + n^{2c_0} L\sqrt k (2L/C_1)^{-k} +Ln^{2c_0}y^{2n} \\
&\ll_{\epsilon,C}& L n^{-c_0}+ (\log n) n^{2c_0} \exp(-(kL)\frac{\log (2L/C_1)}L) + L n^{2c_0}n^{-\frac 2 {C'}}  .
\end{eqnarray*}
We simply choose $L$ large enough (independent of $n$) to ensure that $c_L:=\log (2L/C_1)/L>0$. Note that $kL\ge \epsilon \log n$, therefore 
\begin{eqnarray*}
&& \Prob(\max_{n\le m \le cn} N_m[1-\frac {C\log n}n, 1-\frac{\log n}{C'n}]> \epsilon \log n)  \\
&\ll_{\epsilon,C}& Ln^{-c_0}+ (\log n) n^{2c_0} \exp(-\epsilon c_L \log n) + Ln^{2c_0}n^{-\frac 2 {C'}}  .
\end{eqnarray*}
Thus,  by choosing $c_0<\min(\frac 1 {C'}, \frac{\epsilon c_L}2)$ sufficiently small,  we obtain  
\begin{eqnarray*}
\Prob(\max_{n\le m \le cn} N_m[1-\frac {C\log n}n, 1-\frac{\log n}{C'n}]> \epsilon \log n)  &\ll& n^{-c'}
\end{eqnarray*}
for some $c'>0$.

We now consider $[1-\frac{\log n}{C'n},1]$ where $C'$ may be very large depending on $\epsilon$.  We reset
\begin{eqnarray*}
k&=&1+[\epsilon \log n]
\end{eqnarray*}
and for convenience let $x=1-\frac{\log n}{C'n}$.  Our choice of $C'$ will be large enough to ensure that $C'\epsilon>ce$.  Then
\begin{eqnarray*}
 n^{2c_0}  (\frac{cen(1-x)}{k})^{2k+1} &\le&  n^{2c_0} (\frac{ce}{C'\epsilon})^{2k+1} \\
 &\le & n^{2c_0} n^{-2\epsilon\log(C'\epsilon/(ce))}
\end{eqnarray*}

 Now,  applying   Lemma~\ref{l.near1-max-end} we obtain (for any $c_0>0$ sufficiently small)
\begin{eqnarray*}
&&\Prob(\max_{n\le m\le cn} N_m[1-\frac{\log n}{C'n},1] > \epsilon \log n) \\
&\ll& n^{-c_0} + n^{2c_0} n^{-2\epsilon \log(C'\epsilon/(ce))}+n^{2c_0} n^{-1/C'} .
\end{eqnarray*}
Then,  by choosing $c_0>0$ sufficiently small we  obtain 

\begin{eqnarray*}
\Prob(\max_{n\le m \le cn} N_m[1-\frac{\log n}{C'n},1]> \epsilon \log n)  &\ll& n^{-c'}
\end{eqnarray*}
for some $c'>0$.
 
\subsection{Proof of Lemma~\ref{l.near1-max-sub}} 

Using the fundamental theorem of calculus and Rolle's theorem,  if $p_m$ has at least $k$ zeros in $[x,y]$ then

\begin{eqnarray*}
|p_m(y)|  &\le& \int_{x}^{y} \int_x^{y_1}\dots \int_x^{y_{k-1}} |p_m^{(k)}(y_k)|dy_k\dots dy_1 \\
&\le& I^*_{k}(x,y) 
\ \ :=\ \  \int_{x}^{y} \int_x^{y_1}\dots \int_x^{y_{k-1}} \max_{n\le j\le cn}|p_j^{(k)}(y_k)|dy_k\dots dy_1 
\end{eqnarray*}
therefore
\begin{eqnarray*}
|p_n(y)|  &\le& \max_{n\le m\le cn} |p_m(y)-p_n(y)|+ I^*_{k}(x,y) .
\end{eqnarray*}
We obtain the following estimate:
\begin{eqnarray*}
&& \Prob(\max_{n\le m\le cn} N_m[x,y] \ge \epsilon\log n)   \\
&\le& 
\Prob \Big( |p_n(y)| \ \   \le \ \   I^*_{k}(x,y) + \max_{n\le m\le cn} |p_m(y)-p_n(y)|\Big)
\end{eqnarray*}

 Let $V_n =Var[p_n(y)]=1+y^2 +\dots+y^{2n} \approx (1-y+\frac 1 n)^{-1}$ for $0\le y \le 1$.   When  $1-y \gg 1/n$ we then have $V_n \approx (1-y)^{-1}$, uniformly over $n$.  Thus,  
\begin{eqnarray*}
Var[p_{[cn]}(y)-p_n(y)] &=& (\sum_{n<k\le cn}|y|^{2k}) \\
&\ll&  y^{2n} V_{(c-1)n} \\
&\ll&  y^{2n}V_n.
\end{eqnarray*}
Using Kolmogorov's inequality (see e.g. \cite{Dur}),  it follows that for any $\lambda>0$ we have
\begin{eqnarray*}
\Prob(\max_{n\le m\le cn}  |p_m(y)-p_n(y)|>\lambda) &\ll& \lambda^{-2}Var[p_{[cn]}(y)-p_n(y)] \\
&\ll&  \lambda^{-2} y^{2n}V_n .
\end{eqnarray*}
Thus,
\begin{eqnarray*}
&&\Prob(\max_{n\le m\le cn}  N_m[x,y] \ge \epsilon\log n)   \\
&\le&
\Prob \Big( |p_n(y)|  \le 2\lambda\Big)  + \Prob\Big(I^*_{k}(x,y)>\lambda\Big) + O(\lambda^{-2} y^{2n} V_n).
\end{eqnarray*}
Using Lemma~\ref{l.smallBE},  we have
\begin{eqnarray*}
\Prob(|p_n(y)|\le 2\lambda) & \ll& \frac{\lambda}{\sqrt {V_n}},
\end{eqnarray*}
where we may take  any $\lambda \ge n^{-c_0}$ for some  $c_0>0$ sufficiently small.

By H\"older's inequality and using Doob's $L^2$ maximal inequality (see e.g. \cite{Dur}), we have
\begin{eqnarray*}
&& \mathbb E |I^*_k(x,y)|^2 \\
&\le& \Big(\int_x^y  \dots \int_x^{y_{k-1}} 1dy_k\dots dy_1\Big) \Big(\int_x^y \dots \int_x^{y_{k-1}} \mathbb E [\max_{j\le cn} |p_j^{(k)}(y_k)|^2] dy_k\dots dy_1\Big)\\
&\le& \Big(\frac{(y-x)^k}{k!} \Big) \Big(\int_x^y \dots \int_x^{y_{k-1}} \mathbb E [|p_{[cn]}^{(k)}(y_k)|^2] dy_k\dots dy_1\Big).
\end{eqnarray*}
Therefore using Lemma~\ref{l.iterated-bound} we obtain
\begin{eqnarray*}
 \mathbb E |I^*_k(x,y)|^2
&\ll &   \Big(\frac{(y-x)^k}{k!} \Big) \Big( k! \sqrt{k} (\frac{1}{2(1-y)})^{k+1}\Big) \\  
&\ll&  \frac{ \sqrt{k}}{1-y} (\frac{y-x}{2(1-y)})^{k}  \\
&\ll& \sqrt{k}V_n (2L)^{-k},
\end{eqnarray*}
recalling that   $|y-x|\le |1-y|/L$ by the given assumption.
Using Markov's inequality,  it follows that
\begin{eqnarray*}
\Prob(I^*_k(x,y)>\lambda)  &\ll&  \lambda^{-2}V_n \sqrt{k} (2L)^{-k} 
\end{eqnarray*}
therefore
\begin{eqnarray*}
\Prob(\max_{n\le m\le cn} N_m[x,y] \ge k)  &\ll& \frac{\lambda}{\sqrt{V_n}} + \lambda^{-2} V_n\sqrt{k}(2L)^{-k}+\lambda^{-2} y^{2n} V_n .
\end{eqnarray*}
Choosing $\lambda=n^{-c_0}V_n^{1/2}\ge n^{-c_0}$ for   $c_0 >0$ sufficiently small, we obtain 
\begin{eqnarray*}
\Prob(\max_{n\le m\le cn} N_m[x,y] \ge k) & \ll&  n^{-c_0}  + n^{2c_0}  \sqrt{k}(2L)^{-k}+n^{2c_0}y^{2n}.
\end{eqnarray*}
This completes the proof of Lemma~\ref{l.near1-max-sub}.

\subsection{Proof of Lemma~\ref{l.near1-max-end}} 
We now discuss the proof of the endpoint lemma,  which follows the same steps as the proof of Lemma~\ref{l.near1-max}.  Using the fundamental theorem of calculus and Rolle's theorem,  if $p_m$ has at least $k$ zeros in $[x,1]$ then

\begin{eqnarray*}
|p_m(x)|  &\le& \int_{x}^{1} \int_{x_1}^1\dots \int_{x_{k-1}}^1 |p_m^{(k)}(x_k)|dx_k\dots dx_1 \\
&\le& J^*_{k}(x,1) \ \ 
:= \ \  \int_{x}^{1} \int_{x_1}^1\dots \int_{x_{k-1}}^1  \max_{n\le j\le cn}|p_j^{(k)}(x_k)|dx_k\dots dx_1  
\end{eqnarray*}
therefore
\begin{eqnarray*}
|p_n(x)|  &\le& \max_{n\le m\le cn} |p_m(x)-p_n(x)|+ J^*_{k}(x,1) .
\end{eqnarray*}
We obtain the following estimate:
\begin{eqnarray*}
&&\Prob(\max_{n\le m\le cn} N_m[x,y] \ge \epsilon\log n)  \\
 &\le &
\Prob \Big( |p_n(x)| \ \   \le \ \   J^*_{k}(x,1) + \max_{n\le m\le cn} |p_m(x)-p_n(x)|\Big).
\end{eqnarray*}

Let $V_n =Var[p_n(x)]  \approx (1-x+\frac 1 n)^{-1}$ for $0\le x\le 1$.   Under the given assumption $1\ge 1-x \gg 1/n$ we then have $V_n \approx (1-x)^{-1}$, uniformly over $n$.   Consequently,
\begin{eqnarray*}
Var[p_{[cn]}(x)-p_n(x)] &=& (\sum_{n<k\le cn}|x|^{2k})  \\
&\ll&  x^{2n}V_n.
\end{eqnarray*}

Using Kolmogorov's inequality,  it follows that for any $\lambda>0$ we have
\begin{eqnarray*}
\Prob(\max_{n\le m\le cn}  |p_m(x)-p_n(x)|>\lambda) &\ll& \lambda^{-2}Var[p_{[cn]}(x)-p_n(x)] \\
&\ll&  \lambda^{-2} x^{2n} V_n.
\end{eqnarray*}
Thus,
\begin{eqnarray*}
&&\Prob(\max_{n\le m\le cn}  N_m[x,1] \ge \epsilon\log n)  \\
 & \le &
\Prob \Big( |p_n(x)|  \le 2\lambda\Big)  + \Prob\Big(J^*_{k}(x,x)>\lambda\Big) + O(\lambda^{-2} x^{2n}V_n).
\end{eqnarray*}
Using Lemma~\ref{l.smallBE}, we have
\begin{eqnarray*}
\Prob(|p_n(x)|\le 2\lambda) &\ll& \frac{\lambda}{\sqrt {V_n}},
\end{eqnarray*}
where we may take any $\lambda\ge n^{-c_0}$ for some  $c_0>0$ sufficiently small.

By H\"older's inequality and using Doob's $L^2$-maximal inequality (see e.g. \cite{Dur}), we have
\begin{eqnarray*}
&& \mathbb E[ |J^*_k(x,1)|^2]\\
&\le& \Big(\int_x^1  \dots \int_{x_{k-1}}^1 1dx_k\dots dx_1\Big) \Big(\int_x^1 \dots \int_{x_{k-1}}^1 \mathbb E [\max_{j\le cn} |p_j^{(k)}(x_k)|^2] dx_k\dots dx_1\Big)\\
&\ll& \Big(\frac{(1-x)^k}{k!} \Big) \Big(\int_x^1 \dots \int_{x_{k-1}}^1 \mathbb E [|p_{[cn]}^{(k)}(x_k)|^2] dx_k\dots dx_1\Big).
\end{eqnarray*}
Therefore using Lemma~\ref{l.iterated-bound} we obtain
\begin{eqnarray*}
 \mathbb E |J^*_k(x,1)|^2
&\ll &   \Big(\frac{(1-x)^k}{k!} \Big) \Big((cn)^{2k+1}  \frac{(1-x)^k}{k!}\Big) \\  
&\ll&  \frac{V_n}{(k!)^2} (cn(1-x))^{2k+1}\\
&\ll&  V_n  (\frac{cen(1-x)}{k})^{2k+1}.
\end{eqnarray*}
recalling that $V_n \sim (1-x)^{-1}$. 
Using Markov's inequality,  it follows that
\begin{eqnarray*}
\Prob(J^*_k(x,1)>\lambda) & \ll & \lambda^{-2} V_n  (\frac{cen(1-x)}{k})^{2k+1} 
\end{eqnarray*}
therefore
\begin{eqnarray*}
\Prob(\max_{n\le m\le cn} N_m[x,1] \ge k)  &\ll &  \frac{\lambda}{\sqrt{V_n}} + \lambda^{-2}  V_n  (\frac{cen(1-x)}{k})^{2k+1} +\lambda^{-2} x^{2n}V_n. 
\end{eqnarray*}
Choosing $\lambda=n^{-c_0}V_n^{1/2}\ge n^{-c_0}$ for  $c_0>0$ sufficiently small, we obtain 
\begin{eqnarray*}
\Prob(\max_{n\le m\le cn} N_m[x,1] \ge k)  &\ll & n^{-c_0} + n^{2c_0}  (\frac{cen(1-x)}{k})^{2k+1} +n^{2c_0} x^{2n}.
\end{eqnarray*}
\endproof

\section{Maximal estimates for real roots in $[0,1]$}\label{s.maximal}

In this section,  we will assume that the coefficients $\xi_j$ are independent with zero mean and unit variance and uniformly bounded $(2+\epsilon_0)^{th}$ moments, and their distributions are independent of $n$.

The goal of this section is to prove Lemma~\ref{l.maximal}.  Using Lemma~\ref{l.near0-max} and Lemma~\ref{l.near1-max}, it suffices to show the following estimate.

\begin{lemma} \label{l.maximal-bulk}   
Let  $c>1$ be deterministic constants independent of $n$.  Then for any $C>0$ there are finite constants  $C_0>0$, $C_1>1$, and $C_2>0$, such that 
\begin{eqnarray*}
\Prob\Big(\max_{n\le m\le cn} \Big|N_m [\frac 1 {C_1},1-\frac{C_0\log n}n]-N_{n}[\frac 1 {C_1},1-\frac{C_0\log n}n]\Big| \ge   C_2\Big) 
&\ll& n^{-C} .
\end{eqnarray*}
\end{lemma}

The main ingredient in the proof  is a comparison argument from   \cite{DHV,NNV}, which uses  the following   repulsion estimate  essentially from \cite[Lemma 6]{NNV} and \cite[Theorem 2.2]{DHV}.
\begin{lemma}\label{l.repulsion} For any constant $C>0$, there are constants $C_0>0$, $B>0$, and $C_1>1$, such that  the following holds uniformly over $n\ge 1$:
\begin{eqnarray*}
\Prob(\exists x\in [\frac 1 {C_1},1-\frac{C_0\log n}{n}]: |p_n(x)|\le n^{-B} \text{ and } |p'_n(x)| \le n^{-B}) &\ll& n^{-C}.
\end{eqnarray*}
\end{lemma}
Since this particular formulation of  Lemma~\ref{l.repulsion} was not stated in the above references  for general $\xi_j$ (see \cite[Theorem 2.2]{DHV}),  we sketch a proof of Lemma~\ref{l.repulsion} in the Appendix.

Assuming Lemma~\ref{l.repulsion},  we now prove Lemma~\ref{l.maximal-bulk}.  Let $m$ be any integer such that $n\le m\le cn$.  We shall compare the numbers of   real roots of $p_m$ and $p_n$ inside $I_n:=[\frac 1 {C_1},1-\frac{C_0\log n}{n}]$ by a pairing argument.  The key idea is to show that, with good choice of $C_0,C_1$,  the following holds with overwhelming probablity $1-O(n^{-C})$:
\begin{itemize}
\item $|p_m(x)-p_n(x)|$ is uniformly very small for $x\in I_n$, and 
\item $p_m(x)$ and $p_n(x)$ are unlikely to have double real roots inside $I_n$.
\end{itemize}
Using these properties,  inside $I_n$ we may pair most of the real roots of $p_m$ with most of the real roots of $p_n$.  By combining this argument with an union bound we can ensure that this holds for all $n\le m\le cn$,  still with overwhelming probability $1-O(n^{-C})$.  This argument essentially leads to   the estimates asserted in Lemma~\ref{l.maximal}.

We will use the following comparison lemma. Below, for a real-valued continuous function $f$ and an interval $I$, we denote by $N_f(I)$ the number of real roots of $f$ inside $I$. 
\begin{lemma}\label{l.compare} Let $B$ and $K$ be positive constants.  Then there is a constant $M>0$ that may depend on $B,K$ such that the following holds:

Assume that $I\subset [-1,1]$ and its endpoints may depend on $n\ge 1$.  Assume that $f$ is a real polynomial with $deg(f)\ll n$ and assume that its coefficients are bounded above (in absolute value) by $L \ll n^K$.

Assume that for every $x\in I$ we have $\max(|f(x)|,|f'(x)|)\ge n^{-B}$. 

Then for any real-valued continuous function $g$   on $I$ such that
\begin{eqnarray*}
\sup_{x\in I} |g(x)-f(x)| 
&\ll& n^{-M} 
\end{eqnarray*}
we have
\begin{eqnarray*}
N_f(I) &\le& N_g(I) + O(1).
\end{eqnarray*}
\end{lemma}
\proof
Without loss of generality, we may assume that $n$ and $B$ are sufficiently large,  in particular $B>K+4$.  Then we may take any $M>3B$.

It follows from the given hypothesis that $f$ has no real roots in $I$ with multiplicity $\ge 2$.

We first make several observations.  Let $I'=(a+n^{-B},b-n^{-B})$ where $a<b$ are endpoints of $I$. Then we will show that:
\begin{itemize}
\item[(i)] if $x_1\ne x_2$ are two real roots of $f$ inside $I$ then $|x_1-x_2|> 2n^{-B}$.
\item[(ii)] if $x_1\in I'$ is a root of $f$   then  there is a root $y_1$ of $g$ in $I$ such that $|x_1-y_1| <n^{-2B}$.
\end{itemize}
Using these observations, we may proceed as follows to finish the proof. Let $\eta_1<\dots<\eta_m$ denote the roots of $f$ inside $I$. For each $2\le j \le m-1$,  by (i) we know that $\eta_j \in I'$ and by (ii) there is a root of $g$ inside $(\eta_j - n^{-2B}, \eta_j+n^{-2B}) \subset I$.  Thanks to (i) again, these $n^{-2B}-$neighborhoods of $\eta_j$ are disjoint.    Thus $N_g(I')\ge m-2$ as desired.

We now discuss the proof for (i) and (ii).

For (i),    assume (towards a contradiction) that $f$ has two real roots $x_1<x_2$ inside $I$ with $|x_1-x_2| \le 2n^{-B}$. By the mean value theorem, for some $\alpha$ between $x_1$ and $x_2$ we have $f'(\alpha)=0$. It follows that for every $x\in [x_1,x_2]$ we have
\begin{eqnarray*}
f'(x) &=& f'(\alpha)+O(|x-\alpha|\sup_{y\in I} |f''(y)|) \\
&\ll& |x_1-x_2| n^3L.
\end{eqnarray*}
Since $f'(\alpha)=0$, by the given hypothesis we must have $|f(\alpha)| \ge n^{-B}$. However, by the mean value theorem,
\begin{eqnarray*}
f(\alpha) &=& f(\alpha)-f(x_1) \\
&\ll& |\alpha-x_1| \sup_{x\in (x_1,x_2)} |f'(x)|  \\
&\ll& |x_1-x_2|^2 n^3L\\
&\ll&n^{-2B+3+K}\\
&\ll& n^{-B-1},
\end{eqnarray*}
giving a contradiction.

For (ii), using the given hypothesis we know that $|f'(x_1)|\ge n^{-B}$. Without loss of generality, assume that $f'(x_1) \ge n^{-B}$. Now,  for any $y\in I$ we have $|f''(y)| \ll n^3L \ll n^{2B-1}f'(x_1)$.  It follows that
for $|h|\le n^{-2B}$,
\begin{eqnarray*}
f(x_1+h) &=& f(x_1) + hf'(x_1)+O(h^2\sup_{y\in I}|f''(y)|)\\
&=& hf'(x_1)+n^{-1} O(|h| f'(x_1)).
\end{eqnarray*}

Therefore, for $n$ sufficiently large,
\begin{eqnarray*}
f(x_1-n^{-2B}) &<& -\frac 1 2 n^{-2B}f'(x_1) \ \  < \ \ - \frac 12 n^{-3B} ,\\
f(x_1+n^{-2B}) &>& \frac 1 2 n^{-2B}f'(x_1) \ \ >  \ \ \frac 1 2n^{-3B}.
\end{eqnarray*}
Note that $x_1\pm n^{-2B} \in I$ therefore we may use the comparison condition $\|f-g\|_{sup(I)}=O(n^{-M})=o(n^{-3B})$,
therefore 
\begin{eqnarray*}
g(x_1-n^{-2B}) &<& 0 \ \ < \ \ g(x_1+n^{-2B}),
\end{eqnarray*} and thus the asserted claim will follow by continuity.

\endproof

The following result follows as a corolllary.

\begin{corollary}\label{c.pairing}  Let $B$ and $K$ be positive constants.  Then there is a constant $M>0$ that may depend on $B,K$ such that the following holds:

Assume that $I\subset [-1,1]$ and its endpoints may depend on $n\ge 1$. 
Assume that $f$ and $g$ are real polynomials with degrees bounded above by $O(n)$ and assume that their coefficients are bounded above by $L=O(n^K)$.

Assume that for every $x\in I$ we have $\max(|f(x)|,|f'(x)|)\ge n^{-B}$,  and similarly for every $x\in I$ we have $\max(|g(x)|,|g'(x)|)\ge n^{-B}$.

Assume also that
\begin{eqnarray*}
\sup_{x\in I} |g(x)-f(x)|
& =&O(n^{-M}) .
\end{eqnarray*}

Then
\begin{eqnarray*}
|N_f(I) - N_g(I)| &=& O(1).
\end{eqnarray*}
The implicit constants are independent of $n$ and of $I$.
\end{corollary}

We now discuss the rest of the proof of Lemma~\ref{l.maximal-bulk}.   Fix a constant $C>0$  (which is the desired decay order of the main estimate).

To start, we first truncate $\xi_j$ by a large power of $n$. Let $L=n^K$ where $K>0$ is a sufficiently large constant.  Let 
\begin{eqnarray*}
E &=&\{\max_{n\le j\le cn} |\xi_j| \le L\}.
\end{eqnarray*}
Using the union bound, we have
\begin{eqnarray*}
\Prob(E^c) &=& O(nL^{-2})=O(n^{-(2K-1)}).
\end{eqnarray*}
On the event $E$, we have
\begin{eqnarray*}
\max_{n\le m\le cn} \max_{|x|\le 1}|\sum_{n<j\le m} \xi_j x^{j-n}|   &<& cnL.
\end{eqnarray*}
Therefore, by choosing $K>(C+1)/2$ we obtain
\begin{eqnarray*}
\Prob(\max_{n\le m\le cn}\max_{|x|\le 1}|\sum_{n<j\le m} \xi_j x^{j-n}| \ge cLn) &\le& \Prob(E^c) \\
&\ll& n^{-C}.
\end{eqnarray*}

For $0\le |x|<1-\frac {C_0\log n}n$ we have $x^n=O(n^{-C_0})$, therefore 
\begin{eqnarray*}
(cLn)x^n &\ll&  n^{-C_0+1+K}.
\end{eqnarray*} Consequently,
\begin{eqnarray}\label{e.tailmax}
\Prob(\max_{n\le m\le cn}\max_{|x|\le 1-\frac {C_0\log n}n}|p_m(x)-p_n(x)| \ge n^{-(C_0-K-2)}) &\ll&  n^{-C}.
\end{eqnarray}

On the other hand, using Lemma~\ref{l.repulsion}  and using the union bound,  it follows that there is an event $B_n$ with $P(B_n^c)=O(n^{-C})$ and constants $B>0$, $C_1>1$, and $C_0>0$, such that   on $B_n$ the following holds  for every $m$ integer such that $n\le m\le cn$: 

\begin{eqnarray*}
\max_{x\in  I_m:=[\frac 1 {C_1},  1-\frac{C_0\log m}m]} (m^B |p_m(x)|, m^{B} |p'_m(x)|) &>& 1.
\end{eqnarray*}

Note that $\log x/x$ is a decreasing function (for $x>1$ sufficiently large).   In particular $I_n\subset I_m$ for $m\ge n$ and hence in the above display we may use the same interval $I_n$, for all $n\le m\le cn$.

We now set $M$ to be the constant from Corollary~\ref{c.pairing}, relative to $B$ and $K$.

We certainly have the freedom to choose $C_0$ larger,  thus by making $C_0>K+2+M$ and using \eqref{e.tailmax}, we may assume without loss of generality that on the event $B_n$ we also have,  for every $n\le m\le cn$:
\begin{eqnarray*}
\max_{|x|\le 1-\frac{C_0\log n}{n}} |p_m(x)-p_n(x)| &\le& n^{-M}.
\end{eqnarray*}
Then applying Corollary~\ref{c.pairing} to $f=p_m$ and $g=p_n$, we obtain
\begin{eqnarray*}
\max_{n\le m\le cn}|N_{m}(I_n)-N_{n}(I_n)| &=& O(1).
\end{eqnarray*}

In other words,  for some positive constants $C_0,C_1,C_2$ independent of $n$, we have

\begin{eqnarray*}
\Prob(\max_{n\le m\le cn}|N_{m}[\frac 1 {C_1}, 1-\frac{C_0\log n}{n}]-N_{n}[\frac 1 {C_1}, 1-\frac{C_0\log n}{n}]|>C_2) &\ll& n^{-C}.
\end{eqnarray*}

This completes the proof of Lemma~\ref{l.maximal-bulk}.

\section{Epilogue}\label{s.conclusion}
By symmetry,  the distribution of the real roots of each Kac polynomial is essentially invariant under the symmetry $x\mapsto 1/x$, since $p_n(x)$ and $p^*_n(x):=x^n p_n(1/x)$ have the same distribution.  This well-known observation has  played an important role in many studies for Kac polynomials, reducing many proofs to the fundamental domain $[-1,1]$. However, in the current setting,  this symmetry is no longer applicable. More specifically, while for each $n$ the distribution for $p_n$ is the same as the distribution for $p_n^*$,  in our setting the two sequences $(p_n)_{n\ge 1}$ and  $(p^*_n)_{n\ge 1}$ do not share the same distribution.  In our proof, the pairing argument in our proof of \eqref{e.01} does not seem to extend to the reciprocal polynomial $p^*_n$,  and new ideas are likely needed.  We remark that for the number of real roots on $\mathbb R$, our main result in particular gives an almost sure lower bound
$$\liminf_{n\to\infty} \frac{N_n(\mathbb R)}{\log n}\ge \frac 1 \pi.$$

\appendix
\section{Proof of repulsion estimates}

\subsection{Proof of Lemma~\ref{l.repulsion}} \label{app.a}

In this section,  we sketch the proof of Lemma~\ref{l.repulsion} following the argument in  \cite{DHV, NNV}.  We will assume throughout that $\xi_j$ are independent with zero mean and unit variance and uniformly bounded $(2+\epsilon)^{th}$ moments for some $\epsilon>0$. The implicit constants may depend on $\epsilon$ and the assumed upper bound for the $(2+\epsilon)^{th}$ moments.

For the convenience of the reader, we recall the repulsion estimate.
\begin{lemma}\label{l.repulsion-app} For any constant $C>0$ there are constants $B>0$, $C_0>0$,  and $C_1>1$, such that  the following holds uniformly over $n\ge 1$
\begin{eqnarray*}
\Prob(\exists x\in [\frac 1 {C_1},1-\frac{C_0\log n}{n}]: |p_n(x)|\le n^{-B} \text{ and } |p'_n(x)| \le n^{-B}) 
&\ll& n^{-C}.
\end{eqnarray*}
\end{lemma}

The proof uses the following small ball inequality from \cite[Lemma 7]{NNV}:
\begin{lemma}\label{l.nnv} Let $A>0$.   Then there are constants $C_0>0$ and $C_1>1$ such that the following holds for all $x\in [\frac 1 {C_1},  1-\frac{C_0\log n}n]$ and  $\delta=n^{-A}$:
\begin{eqnarray*}
\Prob(|p_n(x)|\le \gamma) &\ll& \gamma^{0.99}.
\end{eqnarray*}
\end{lemma}
 
Now, let   $K>\frac {C+1}2$ be a large constant and define
\begin{eqnarray*}
E &=& \{\max_{0\le j\le n}|\xi_j|\le n^K\}.
\end{eqnarray*}
Using the union bound, we have $\Prob(E^c)=O(nn^{-2K})=O(n^{-C})$.

We will then choose $A>0$ sufficiently large, in particular $A>100(K+3)$ and $A>2C$. 

Let $B>2A$ and let $C_0, C_1$ be given by Lemma~\ref{l.nnv} using $\delta'=n^{-1.99A}$.

We now divide $[\frac 1 {C_1},1-\frac{C_0\log n}{n}]$ into subintervals of length $\delta=n^{-A}$ where $A>100(K+3)$.  Such division will need at most $O(\delta^{-1})$ intervals.  If $I$ is one such interval (whose center is denoted by $x_I$) and if there is some $x\in I$ with 
\begin{eqnarray*}
\max(|p_n(x)|, |p'_n(x)|)  &\le& n^{-B}
\end{eqnarray*}
 then conditioning on $E$ we will show that
\begin{eqnarray*}
|p_n(x_I)| \ll  \delta^{1.99}.
\end{eqnarray*}
Indeed,  we may use Taylor expansion around $x$ to obtain 
\begin{eqnarray*}
|p_n(x_I)| &\le& |p_n(x)|+|(x_I-x)p_n'(x)| + \frac 1 2 |x_I-x|^2 \sup_{t\in I}| p_n''(t)|  \\
&\le& n^{-B} + \frac  1 2 \delta n^{-B} + \frac 1 2 \delta^2 (n^{3+K}) \\
&\le& \delta^2 + \frac  1 2 \delta^3   + \frac 1 2 \delta^2 \delta^{-1/100}  \\
&\ll& \delta^{1.99}.
\end{eqnarray*}

Consequently, using the union bound, we obtain
\begin{eqnarray*}
&&\Prob(\exists x\in [\frac 1 {C_1},1-\frac{C_0\log n}{n}]: |p_n(x)|\le n^{-B} \text{ and } |p'_n(x)| \le n^{-B})\\
&\le& \Prob(E^c) +\sum_{I} \Prob(|p_n(x_I)|\ll \delta^{1.99}) \\
&\ll& n^{-C} + \delta^{-1}(\delta^{1.99})^{0.99}  \ \ \ \text{(by Lemma~\ref{l.nnv})}\\
&\ll& n^{-C}.
\end{eqnarray*}
This completes the proof of Lemma~\ref{l.repulsion}.

\end{document}